\documentclass[10pt]{article}
\usepackage{amsmath} \usepackage{amsfonts} \usepackage{amssymb}
\usepackage{textcomp}
\usepackage{graphicx}
\usepackage{color}

\DeclareMathSymbol{\leqslant}{\mathalpha}{AMSa}{"36} 
\DeclareMathSymbol{\geqslant}{\mathalpha}{AMSa}{"3E} 
\DeclareMathSymbol{\eset}{\mathalpha}{AMSb}{"3F}     
\renewcommand{\leq}{\;\leqslant\;}                   
\renewcommand{\geq}{\;\geqslant\;}                   

\makeatletter
\def\captionfont@{\footnotesize}
\def\captionheadfont@{\scshape}
\long\def\@makecaption#1#2{%
  \vspace{2mm}
  \setbox\@tempboxa\vbox{\color@setgroup
    \advance\hsize-6pc\noindent
    \captionfont@\captionheadfont@#1\@xp\@ifnotempty\@xp
        {\@cdr#2\@nil}{.\captionfont@\upshape\enspace#2}%
    \unskip\kern-6pc\par
    \global\setbox\@ne\lastbox\color@endgroup}%
  \ifhbox\@ne 
    \setbox\@ne\hbox{\unhbox\@ne\unskip\unskip\unpenalty\unkern}%
  \fi
  \ifdim\wd\@tempboxa=\z@ 
    \setbox\@ne\hbox to\columnwidth{\hss\kern-6pc\box\@ne\hss}%
  \else 
    \setbox\@ne\vbox{\unvbox\@tempboxa\parskip\z@skip
        \noindent\unhbox\@ne\advance\hsize-6pc\par}%
\fi
  \ifnum\@tempcnta<64 
    \addvspace\abovecaptionskip
    \moveright 3pc\box\@ne
  \else 
    \moveright 3pc\box\@ne
\nobreak
\vskip\belowcaptionskip
\fi
\relax
}
\makeatother

\def\writefig#1 #2 #3 {\rlap{\kern #1 truecm
\raise #2 truecm \hbox{#3}}}

\newtheorem{thm}{Theorem}[section]
\newtheorem{lem}[thm]{Lemma}

\newtheorem{rem}[thm]{Remark}

\newcommand{\bbC}{{\ensuremath{\mathbb C}} }
\newcommand{\bbD}{{\ensuremath{\mathbb D}} }
\newcommand{\bbE}{{\ensuremath{\mathbb E}} }
\newcommand{\bbF}{{\ensuremath{\mathbb F}} }

\newcommand{\bbQ}{{\ensuremath{\mathbb Q}} }

\newcommand{\bbU}{{\ensuremath{\mathbb U}} }
\newcommand{\bbV}{{\ensuremath{\mathbb V}} }

\def\1{\ifmmode {1\hskip -3pt \rm{I}} \else {\hbox {$1\hskip -3pt \rm{I}$}}\fi} 

\newcommand{\df}{\stackrel{\Delta}{=}}

\newcommand{\case}[1]{C{\small ASE}~#1}

\newcommand{\abs}[1]{\lvert#1\rvert}  

\newcommand{\setof}[2]{\left\{#1\,:\,#2\right\}}

\newcommand{\bracket}[1]{\left\{ #1\right\}}
\newcommand{\dif}{{\rm d}}

\begin{document}
\title{Entropy-driven phase transition in a polydisperse hard-rods
lattice system}
\author{
D. Ioffe\\
Faculty of Industrial Engineering, Technion\\
Haifa 32000, Israel\\
\texttt{ieioffe@ie.technion.ac.il}
\and
Y.Velenik\\UMR-CNRS 6085, Universit\'e de Rouen\\F-76821 Mont Saint Aignan, France\\
\texttt{Yvan.Velenik@univ-rouen.fr}
\and
M. Zahradn\'ik\\
Faculty of Mathematics and Physics, Charles University\\
18600 Praha 8, Czech Republic\\
\texttt{mzahrad@karlin.mff.cuni.cz}
}
\maketitle

\begin{abstract}
We study a system of rods on $\mathbb{Z}^2$, with hard-core exclusion. Each rod
has a length between $2$ and $N$. We show that,
when $N$ is sufficiently large, and for suitable fugacity, there are several
distinct Gibbs states, with orientational long-range order. This is in sharp
contrast with the case $N=2$ (the monomer-dimer model), for which Heilmann and
Lieb proved absence of phase transition at any fugacity.
This is the first example of a pure hard-core system with phases displaying
orientational order, but not translational order; this is a
fundamental characteristic feature of liquid crystals.
\end{abstract}

\section{Introduction and results}
In 1949, Lars Onsager proposed a theory of the isotropic-nematic phase
transition in liquid crystals, which relied on the following simple
heuristics~\cite{On49}. Picture each molecule as a (very) long, (very)
thin rod.
There is no energetical interaction between the rods, except for hard-core
exclusion. Since at low densities the molecules are typically far from each
other, the resulting state will be an isotropic gas. However, at large
densities
it might be more favorable for the molecules to align spontaneously, since the
resulting loss of orientational entropy is by far compensated by the gain of
translational entropy: indeed, there are many more ways of placing nearly
aligned rods than randomly oriented ones.

This is probably the first example of an entropy-driven phase transition. It
shows that an increase of entropy can sometimes result in an apparently more
ordered structure, hence the often used expression ``order from disorder''.

\medskip
In spite of the obvious physical relevance of such issues, rigorous results are
still very scarce. The only proof of such a phase transition has been given
(as a side remark) in~\cite{BrKuLe84} for the following simple model: The rods
are one-dimensional unit-length line segments in $\mathbb{R}^2$, with two
possible orientations (say, horizontal and vertical); a configuration of $N$
rods ($N$ is not fixed) is specified by a family
 $(x,\sigma) \in \mathbb{R}^{2N}
\times \{-1,1\}^N$, where $(x_{2k-1}, x_{2k})$ is the position of the middle of
the $k$th rod, while $\sigma_k$ represents its orientation. A configuration in a
subset $V$ of $\mathbb{R}^2$ is admissible if all rods are inside $V$ and are
disjoint; let $\mathcal{C}_V$ denote this event. Then the measure describing
the process is $\nu_\lambda(\,\cdot\,|\, \mathcal{C}_V)$, where $\nu_\lambda$ is
the product of the Poisson point process in $\mathbb{R}^2$ of intensity
$\lambda>0$ and the Bernoulli process of parameter $\tfrac12$. The main result
is then that, in the thermodynamic limit $V\nearrow \mathbb{R}^2$, there are (at
least) two limiting Gibbs states with long-range orientational order, for all
$\lambda$ large enough. This model has a very special
feature though, namely two horizontal (respectively. vertical) rods have 0-probability
of intersecting under $\nu_\lambda$. This is a considerable simplification, and
therefore this result does not provide any information on the much more
interesting case of rods of finite width.

\medskip
The other rigorous results are concerned with lattice versions of this
problem.
Heilmann and Lieb proved in a classical paper~\cite{HeLi72} that there
is \textit{no} phase transition in the monomer-dimer model, in the sense that
the corresponding free energy is always analytic. This model is defined as
follows (on $\mathbb{Z}^2$, their result holds more generally). Let $V$ be a
finite subset of $\mathbb{Z}^2$. For $x,y\in\mathbb{Z}^2$, we write $x\sim
y$ if $|x-y|=1$. Let $\mathcal{E_V} = \{\{x,y\}\subset V \,:\, x\sim
y\}$; for $e,e'\in\mathcal{E_V}$ we write $e\sim e'$ if $e\cap e'\neq\emptyset$.
The state space is given by $\Omega= \{0,1\}^\mathcal{E_V}$. A configuration
$\omega\in\Omega$ is admissible if $e\sim e' \implies
\min(\omega_e,\omega_{e'})=0$. The probability of a configuration $\omega$ is
then given by
$$
\mu_\lambda(\omega) \propto \mathbf{1}_{\{ \omega \text{ is admissible} \}} \;
\lambda^{|\omega|}\,,
$$
where $|\omega| = \sum_{e\in\mathcal{E}_V} \omega_e$, and $\lambda>0$.
Informally, when a pair of sites $e$ is such that $\omega_e = 1$, then the two
sites are occupied by a dimer; a configuration is admissible if no site
belongs to more than one dimer; $\lambda$ is the dimers' fugacity.

An alternative approach to this model was then discovered by van
den Berg~\cite{vdB99}. Using disagreement percolation methods, he was able to
give a very simple proof of the much stronger result that this model is in fact
completely analytic, in the sense of~\cite{DoWa90}. This paper is also
interesting in that it clearly points out the very special nature of dimers.
Indeed, it would be impossible to push the analysis to arbitrary values of
fugacities, were it not for a magical property of dimers: A monomer-dimer model
on a graph $G$ is actually equivalent to a pure hard-core gas on the line-graph
of $G$.

\medskip
Several other similar models have been introduced (see, e.g.,
\cite{LeGa71,He72,Hu77}), in which existence of orientationally ordered states
has been proven. All these models, however, share the same defect, namely the
ordered states also automatically display long-range translational order, i.e.
they are perturbations of periodic configurations. Thus they really can't be
considered satisfactory models of liquid crystals, since a central characteristic
of the latter is the liquid-like spatial behavior in the ordered states. In
order to solve this problem, Heilmann and Lieb~\cite{HeLi79} proposed five
different models of hard-core particles (actually dimers). For these models they
proved the existence of long-range orientational order at low temperatures, and
gave quite plausible arguments in favour of the absence of long-range
translational order. These, however, were not pure hard-core
models, since an additional attractive interaction favouring alignment of the
dimers was introduced, and thus the question of whether pure hard-core
interaction can give rise to such phases was left open (actually, Heilmann and
Lieb even stated that it was ``\textsl{doubtful [...] whether hard rods on a
cubic lattice without any additional interaction do indeed undergo a phase
transition}'').

\medskip
To the best of our knowledge, these are the only rigorous results pertaining to
this problem. It would be extremely desirable to prove the
existence of a phase transition in the monomer--$k$-mer model (replacing dimers
above by $k$-mers, i.e. families of $k$ aligned nearest-neighbor sites), for
large enough $k$. This seems rather delicate however, and in this work we
concentrate on another variant of the monomer-dimer model, with only hard-core
exclusion and for which it is actually possible to prove existence of phases
with orientational long-range order and no translational long-range order;
actually it is also seems possible to treat the three-dimensional case, which 
presumably would lead to 
a counterexample to the above claim.  We hope to return to the
monomer--$k$-mer problem and to the case of higher dimensions 
in the future.

\bigskip
Our model is defined as follows. We call rod a family
of $k$, $k\in\mathbb{N}$, distinct, aligned, nearest-neighbor sites of
$\mathbb{Z}^2$ a $k$-rod is a rod of length $k$, and we refer to $1$-rod as
vacancies. Let $V\subset\mathbb{Z}^2$; a configuration $\omega$ of our model
inside $V$ is a partition of $V$ into a family of disjoint rods. We write
$N_k(\omega)$ the number of $k$-rods in $\omega$. The probability of the
configuration $\omega$ is given by
\begin{equation}
\label{mu_measure}
\mu_{q,N,V}(\omega) \propto
\mathbf{1}_{\{N_k(\omega)=0, \, \forall k>N\}}(2q)^{N_1(\omega)}\,
q^{\sum_{k=2}^N N_k(\omega)}\,,
\end{equation}
where $q>0$ and $N\in\mathbb{N}$. Informally, only rods of length at most $N$
are allowed; the activity of each rod of length at least $2$ is $q$, and is
independent of the rod's length; there is an additional activity $2q$ for
vacancies.

\begin{rem}
The activity of vacancies can be removed at the cost of introducing an
additional factor $(2q)^{-k}$ for each $k$-rods ($k\geq 2$).
\end{rem}

Our main task is the proof of the following theorem, which states that for large
enough $N$, there is a phase transition from a unique (necessarily isotropic)
Gibbs state at large values of $q$ to several Gibbs states with long-range
orientational order at small values of $q$, but \textit{no} translational order.
This is thus the first model, where such a behavior can be proved.
\begin{thm}
\label{thm_Main}
\begin{enumerate}
\item For any $N\geq 2$, there exists 
$q_0=q_0(N)>0$ such that, for all
$q\geq q_0$ there is a unique, isotropic Gibbs state.
\item 
For any $q >0$ sufficiently small
there exist $N_0 = N_0 (q)$, 
such 
that for all $N\geq N_0$ 
there are two different extremal Gibbs states with
long-range orientational order. More precisely, there exists a Gibbs state
$\mu_{q,N}^{h}$ such that
\begin{equation}
\label{mu_horizontal}
\mu_{q,N}^{h}(\,0 \text{ belongs to a horizontal rod }\,) > 1/2\,.
\end{equation}
\end{enumerate}
\end{thm}

In the sequel we shall refer to the infinite volume Gibbs state 
$\mu_{q,N}^{h}$ as to the horizontal state. By symmetry the $\pi /2$ rotation
of the latter gives the vertical Gibbs state 
$\mu_{q,N}^{v}$, which would statistically favour vertically oriented  
rods. 

A funny consequence of the techniques we develop in order to prove
 Theorem~\ref{thm_Main} is the following result on a sampling of
infinite volume horizontal and vertical states by the shapes of the
family of finite volume domains:

\begin{thm}
\label{thm_Domains}
Let $\underline{k}= (k_1, k_2)$ be two 
natural numbers. For $n=1,2,\dots$, consider lattice 
boxes
\[
V_n^{\underline{k}} = [-k_1n,\dots ,k_1 n]\times [-k_2n,\dots ,k_2 n] ,
\]
and let $\mu_{q,N ,V_n^{\underline{k}}}$ be the finite volume  Gibbs state 
specified in \eqref{mu_measure}. Then if $q$ and $N$ satisfy conditions
of 2) of Theorem~\ref{thm_Main}, 
\begin{equation*}
\label{Domains}
\begin{split}
\lim_{n\to\infty}  \mu_{q,N ,V_n^{\underline{k}}} \, &=\,
\mu_{q, N}^{h}\quad \text{if}\ k_1 > k_2\\
&\text{and}\\
\lim_{n\to\infty}  \mu_{q,N ,V_n^{\underline{k}}} \, &=\,
\mu_{q, N}^{v}\quad \text{if}\ k_1 < k_2
\end{split}
\end{equation*}
\end{thm}

Theorem~\ref{thm_Main} is proved by showing that, for $N$ large enough, the
model defined above is a small perturbation (in a suitable sense) of the
``exactly solvable'' case $N=\infty$. For the latter, the theorem takes the
following form. Let $q_{\textrm{c}} = 1/(2+2\sqrt{2})$.
\begin{thm}
\label{thm_exact}
\begin{enumerate}
\item Let $N=\infty$. For all $q\geq q_{\textrm{c}}$ there 
is a unique, isotropic
Gibbs state.
\item For all $q<q_{\textrm{c}}$, there are (at least) $2$ different extremal
Gibbs states with long-range orientational order. More precisely, 
there exists a
Gibbs state $\mu_q$ such that
$$
\mu_q(\,0 \text{ belongs to a horizontal rod }\,) > 1/2\,.
$$
\end{enumerate}
\end{thm}

\section{An exactly solvable case: Proof of Theorem~\ref{thm_exact}}
\label{ssec_exact}
In this section, we show that the model obtained by setting $N=\infty$ is
actually exactly solvable, since it can be mapped on the 2D Ising model.

We suppose that our system is contained inside a square box $V$ of linear size
$L$; we suppose that we have periodic boundary conditions. We want to partition
$V$ into two disjoint subsets corresponding to the regions occupied by
horizontal and vertical rods respectively. This can be done easily once we have
said what we do with vacancies. The trick is to split vacancies into two
species, horizontal and vertical. Doing so, starting from a configuration
$\omega$ of our original model, we obtain a family of $2^{N_1(\omega)}$
different configurations $\widetilde{\omega}_i$, $i=1,\ldots,N_1(\omega)$. The
probability of each such configuration is then taken to be
\begin{equation*}
\mu_{q,N,V}(\widetilde{\omega}) \propto
\mathbf{1}_{\{N_k(\omega)=0, \, \forall k>N\}}\,
2^{-N_1(\widetilde\omega)}\, (2q)^{N_1(\widetilde\omega)}\, q^{\sum_{k=2}^\infty
N_k(\widetilde\omega)}
= q^{\sum_{k=1}^\infty N_k(\widetilde\omega)}\,.
\end{equation*}
We can now partition $V = V_{\mathrm{h}} \vee V_{\mathrm{v}}$ into two disjoint
subsets. Once these subsets are fixed, the problem is reduced to the study of
one-dimensional partition functions; indeed each maximal connected
horizontal piece of $V_{\mathrm{h}}$ can be filled by horizontal rods
independently of what choice is made for the rest of the configuration, and
similarly for vertical pieces of $V_{\mathrm{v}}$.

In the $N=\infty$ the one-dimensional partition functions 
 could be 
computed exactly: 
\begin{equation*}
Z^{\mathrm{1D}}_n
= \sum_{i=0}^{n-1} \binom{n-1}{i}\, q^{i+1}
= q\, \left( 1 + q \right)^{n-1}
= \frac q{ 1 + q }\, \left( 1 + q \right)^{n}\,,
\end{equation*}
where $Z^{\mathrm{1D}}_n$ is the 1D partition function in a box of length
$n\geq 1$.
Now observe that the exponentially decreasing term $\left( 1 + q \right)^{n}$ is
actually irrelevant, since its total contribution to the weight of a partition
$V_{\mathrm{h}} \vee V_{\mathrm{v}}$ is
$
\left( 1 + q \right)^{|V|}\,,
$
and is therefore independent of the partition. We thus see that this model
possess the remarkable property that all its 1D partition functions are actually
equal to $e^{-4\beta} \stackrel{\triangle}{=} q / ( 1 + q )$, and
therefore independent of $n$. It is then very easy to compute the total weight
of a partition:
\begin{equation*}
\mathrm{weight}(V_{\mathrm{h}},V_{\mathrm{v}}) \propto
e^{-2\beta|\underline\gamma|}\,,
\end{equation*}
where $\underline\gamma=(\gamma_1,\ldots,\gamma_m)$ is the set of contours of
the partition, i.e. the set of all bonds of the dual lattice intersecting a bond
between two nearest-neighbor sites belonging one to $V_{\mathrm{h}}$ and the
other to $V_{\mathrm{v}}$; $|\underline\gamma|$ is the total length of the
contours, where the two components of the partition are reinterpreted as the
components occupied by $+$, respectively. $-$ spins.

One thus observes that the weight of partitions are the same as those of the
corresponding configuration of the 2D Ising model at inverse temperature
$\beta$, in the box $V$ with periodic boundary conditions. Now notice that
$\beta(q_{\textrm{c}}) = \tfrac12 \log(1+\sqrt{2})= \beta_{\textrm{c}}$, the
critical inverse temperature of the 2D Ising model. It immediately follows that
for $q \geq q_{\textrm{c}}$, the corresponding Ising model is in the
high-temperature phase, and therefore possesses a unique Gibbs state. Statement
1 of Theorem~\ref{thm_exact} follows immediately from the symmetry of the
latter Gibbs state.

To prove statement 2. requires only a simple additional argument. For a given
collection of rods $\widetilde\omega$, let us denote by $Z_{\mathrm{h}}$ the
number of sites of $V_{\mathrm{h}}$ containing vacancies, and $N_h =
|V_{\mathrm{h}}| - Z_{\mathrm{h}}$; similarly introduce $N_{\mathrm{v}}$ and
$Z_{\textrm{v}}$. When $q<q_{\textrm{c}}$, the Ising model is in the
low-temperature region; consequently,
\begin{equation*}
\mathbb{E}_{\mu_{q,V}} \bigl[ \bigl\vert |V_{\mathrm{h}}| - |V_{\mathrm{v}}|
\bigr\vert \bigr]
=
\mathbb{E}_{\textrm{Ising},\beta,V} \bigl[ \bigl\vert \sum_{x\in V} \sigma_x
\bigr\vert \bigr]
> c L^2\,,
\end{equation*}
with $c>0$. Since
$
\mu_{q,V}(\,0 \text{ belongs to a horizontal rod }\,) =
L^{-2}\, \mathbb{E}_{\mu_{q,V}} \bigl[ N_{\mathrm{h}} \bigr]\,,
$
the conclusion now follows easily from
$$
\bigl\vert |V_{\mathrm{h}}| - |V_{\mathrm{v}}| \bigr\vert \leq
\bigl\vert N_{\mathrm{h}} - N_{\mathrm{v}} \bigr\vert +
\bigl\vert Z_{\mathrm{h}} - Z_{\mathrm{v}} \bigr\vert\,,
$$
and $\mathbb{E}_{\mu_{q,V}} \bigl[ \bigl\vert Z_{\mathrm{h}} -
Z_{\mathrm{v}} \bigr\vert \bigr] < C L$, by the Central Limit Theorem.

\begin{rem}
A lot of additional information (e.g., on the critical behavior) can be
extracted from this mapping to the 2D Ising model. We refrain from doing that
here, since this is quite straightforward...
\end{rem}
\begin{rem}
As it was pointed to one of us by Lincoln Chayes, in $N=\infty$ case 
the techniques of reflection positivity enable to treat a more 
general situation when the rod weights are 
given by
\[
\lambda^{N_1 (\omega  )}\prod_k  q^{N_k (\omega )} .
\]
In the above notation the case we consider here corresponds to a 
specific choice $\lambda = 2q$. However, the reflection positivity
argument does not go through when there is a finite collection
of admissible rod lengths, $N<\infty$. 
\end{rem}

\section{Asymptotics of one-dimensional partition functions}
\label{asymptotics1DPF}
Our next step is to show that the model with finite (but large) $N$ is actually
a small perturbation of the exactly solvable model analyzed in
Section~\ref{ssec_exact}. The idea, which is described in details in
Subsection~\ref{ssec_PS} is to replace all the 1D partition functions by their
limiting values (for $n\to\infty$), and to expand the error term. To be able to
control this expansion, we need a very good control on the speed of convergence
of these 1D partition functions. This is the aim of the current subsection.

\subsection{The setup.}
We shall consider here a general case of non-negative
rod activities $\{f_k\}$ which we shall view as a perturbation
of the geometric distribution,
\begin{equation}
\label{f_type}
f_k = qp^{k-1} +\epsilon_k ;\ \ k=1,2,\dots ,
\end{equation}
where $p+q =1$ and the activities $\{f_k\}$  are normalized to
furnish a probability distribution, that is
\begin{equation}
\label{eps_sum}
\sum_k \epsilon_k = 0
\end{equation}
 The important assumptions are those on the smallness of the
perturbation sequence $\{\epsilon_k\}$ with respect to the
background geometric distribution $\{ qp^{k-1}\}$:

\medskip
\noindent
{\bf Assumption A1} There exist $\delta <\infty$ and
$\rho\in (1, p^{-1}]$ such that
\begin{equation}
\label{eps_bound}
|\epsilon_k |\leq \delta \rho^{-k} ,\ \ k=1,2,\dots .
\end{equation}
 \medskip
\noindent
{\bf Assumption A2} There exists $\alpha >0$ sufficiently small
such that,
\[
\delta < \alpha (\rho -1 )^2 .
\]

\noindent
Assumption~{\bf A1} is an essential one. On the other hand,
Assumption~{\bf A2} is more technical and it merely
reflects an intended
compromise between giving a relatively simple proof and
yet generating a whole family of examples where the entropy driven
phase transition takes place. Notice that since $\rho <1/p$, assumption
 {\bf A2} in fact implies a bound on $\delta$ in terms of $q$:
\[
\delta < \frac{\alpha}{p^2}q^2 .
\]

Given $\bracket{ f_k}$ as in \eqref{f_type} above we use it
to set up the renewal relation:
\begin{equation}
\label{renewal}
g_0 =1\qquad {\rm and}\qquad g_n =\sum_{k=1}^n f_k g_{n-k}
\qquad \ n=2,3,\dots .
\end{equation}
Define the generating function of the $\bracket{f_k}$ sequence as
\[
{\mathbb F}(\xi) = \sum_{k=1}^\infty f_k \xi^k  = \frac{q\xi}{1-p\xi} +
\sum_k \epsilon_k\xi^k \df \bbQ (\xi )+\bbE (\xi ).
\]
Above $\bbQ $ is the generating function of the geometric
distribution $\{ qp^{k-1}\}$ and, accordingly, $\bbE$ is the generating
function of $\{\epsilon_k\}$.

In the sequel we use the notation
\[
{\mathbb D}_r (x) = \setof{z\in \bbC}{|z-x|<r}
\]
for an open complex ball of radius $r$ centered at $x$.  By
 {\bf A1}, ${\mathbb F}$ is analytic on $\bbD_\rho (0)$.

 The generating function ${\mathbb G}$
of the $\bracket{g_n}$ sequence is defined and analytic in $\setof{z\in
{\mathbb C}}{\abs{z}<1}$.
 By
the usual renewal theory,
\begin{equation}
\label{FN_limit}
\lim_{n\to\infty} g_n \, =\, \frac1 {{\mathbb F}_N^{\prime} (1)}\,=
\frac1{1/q +\sum_k k\epsilon_k}\,
\df\, g .
\end{equation}
\subsection{The representation formula}
For every $\nu \in (0,1)$,
\begin{equation}
\label{represent}
\begin{split}
g_n -g \, &=\, \frac1{2\pi i}\oint_{D_\nu (0)}
\bracket{
\frac{\dif \xi}{\xi^{n}(1-{\mathbb F} (\xi ))} -
\frac{\dif \xi}{\xi^{n} (1-\xi ) {\mathbb F}^{\prime} (1)}}\\[2ex]
&= \frac1{2\pi i}\oint_{D_\nu}\frac{\dif \xi }{\xi^n}
\bracket{\frac{({\mathbb F} (\xi ) -1) - {\mathbb F}^{\prime} (1)
(\xi -1)}{({\mathbb F} (\xi ) -1)(\xi -1) {\mathbb F}^{\prime} (1)}}\\
\end{split}
\end{equation}
Now,
\[
({\mathbb F} (\xi ) -1) - {\mathbb F}^{\prime} (1) =
\left\{({\mathbb Q} (\xi ) -1) - {\mathbb Q}^{\prime} (1)\right\}
+
\left\{ \bbE (\xi ) -\bbE^\prime (1) (\xi -1)\right\}
\]
As a result,
\begin{equation}
\label{eq:U1}
\begin{split}
2\pi i \bbF^{\prime} (1)\left(
 g_n -g\right)\, &=\,
\oint_{D_\nu}\frac{\dif \xi }{\xi^n}
\bracket{\frac{({\mathbb Q} (\xi ) -1) - {\mathbb Q}^{\prime} (1)
(\xi -1)}{({\mathbb F} (\xi ) -1)(\xi -1) }}
+
\oint_{D_\nu}\frac{\dif \xi }{\xi^n}
\bracket{ \frac{\bbE (\xi ) - \bbE^{\prime} (1)(\xi -1)}{(\bbF (\xi )-1)
(\xi -1)}} \\[2ex]
&\df
\oint_{D_\nu}\frac{\dif \xi }{\xi^n}
\bracket{\frac{({\mathbb Q} (\xi ) -1) - {\mathbb Q}^{\prime} (1)
(\xi -1)}{({\mathbb F} (\xi ) -1)(\xi -1) }}
+
\oint_{D_\nu}\frac{\dif \xi }{\xi^n}\bbU_1 (\xi ) .
\end{split}
\end{equation}
On the other hand, since the geometric distribution $\{ qp^{k-1}\}$
generates (via the renewal relation) a constant sequence $\{ q\}$,
\[
\oint_{D_\nu}\frac{\dif \xi }{\xi^n}
\bracket{\frac{({\mathbb Q} (\xi ) -1) - {\mathbb Q}^{\prime} (1)
(\xi -1)}{({\mathbb Q} (\xi ) -1)(\xi -1) }} \equiv 0 .
\]
Subtracting the above expression for zero from the first term
on the right hand side of \eqref{eq:U1} we obtain
\begin{equation}
\label{eq:U2}
 -
\oint_{D_\nu}\frac{\dif \xi }{\xi^n}
\bracket{
\frac{\left({\mathbb Q} (\xi ) -1- {\mathbb Q}^{\prime} (1)
(\xi -1)\right)\bbE (\xi )}{({\mathbb F} (\xi ) -1)
({\mathbb Q} (\xi ) -1)(\xi -1) }} \df
\oint_{D_\nu}\frac{\dif \xi }{\xi^n}\bbU_2 (\xi ) .
\end{equation}
Thus, with $\bbU_1$ and $\bbU_2$ being defined as in
\eqref{eq:U1} and \eqref{eq:U2} above, the representation
formula \eqref{represent} reads as
\[
2\pi i \bbF^{\prime} (1)\left(
 g_n -g\right)\, =\,
\oint_{D_\nu}\frac{\dif \xi }{\xi^n}
\bracket{ \bbU_1 (\xi ) +\bbU_2 (\xi )} .
\]

Assume that it so happens that both ${\mathbb U}_1 $ and
$ \bbU_2$ are analytic in an
open neighbourhood of
$\bbD_R (0)$ for some $R>1 $.
 Then \eqref{represent} implies that
\begin{equation}
\label{exp_bound}
\abs{ g_n - g}\, \leq\, \left( \max_{\abs{\xi}=R}\,
\abs{{\mathbb U_1} (\xi )}+\max_{\abs{\xi}=R}\,
\abs{{\mathbb U_2} (\xi )}\right)\, \frac1 {2\pi\bbF^\prime (1) R^n} .
\end{equation}
We shall represent $\bbU_1$ and $\bbU_2$ as ratios of two analytic functions,
In Subsection~\ref{sub:lower} we derive a lower bounds on the
denominators, whereas in Subsection~\ref{sub:upper} we derive the
corresponding upper bound for the numerators. Eventually we shall
pick $R= (1+\rho )/2$ and the target bound on
 $\left(\max_{\abs{\xi}=R}
\abs{{\mathbb U} (z)}+\max_{\abs{\xi}=R}\abs{{\mathbb U_2} (\xi)}\right)
$  is formulated in Subsection~\ref{sub:target}.
Finally the case of uniform rod weights is worked out in detail in
Subsection~\ref{sub:uniform}

\subsection{Lower bounds on
the denominators}
\label{sub:lower}
By {\bf A1} the function
\[
\frac{\bbF (\xi )- 1}{\xi -1} = \frac{1}{1- p\xi} + \frac{\bbE (\xi )}{\xi-1}
\df \frac{1}{1- p\xi} +\bbV (\xi )
\]
is analytic in $\bbD_\rho (0)$.

Let us pick $\eta \in (0,\rho -1 )$, later on we shall settle down
with the choice $\eta =(\rho -1)/2$, but in principle all the estimates
below could be further optimized. Since $1+\eta <\rho < 1/p$,
\begin{equation}
\label{lower}
\inf_{\xi\in\bbD_{1+\eta} (0) }\left|  \frac{1}{1- p\xi}\right| \geq
\frac{1}{1+(1+\eta) p} \geq \frac12 .
\end{equation}
It, therefore, remains to derive an appropriate upper bound on
$|\bbE (\xi )/(\xi-1)| = |\bbV (\xi )|$. There are two cases to be considered:

\medskip
\noindent
\case{1}. $\xi\in \bbD_{1+\eta}(0)\setminus\bbD_\eta (1)$. Then, by
{\bf A1},
\begin{equation}
\label{upper1}
\left|
\bbV (\xi )\right| \leq \frac{\delta}{\eta}
\sum_1^\infty \left( \frac{1+\eta}{\rho}\right)^k =
\frac{\delta(1+\eta)}{\eta (\rho -(1+\eta ))}  .
\end{equation}

\medskip
\noindent
\case{2}. $\xi\in \bbD_{\eta}(1)$. Since $\bbE (\cdot )$ is analytic
in $\bbD_{\eta}(1)$,
\[
\bbE (\xi ) = \sum_1^\infty\epsilon_k\xi^k = \sum_1^\infty\epsilon_k
\sum_{l=0}^k {k\choose l} (\xi -1)^l =
\sum_1^\infty \tilde{\epsilon}_l
 (\xi -1)^l ,
\]
 where we have used $\sum\epsilon_k =0$ and, accordingly, have defined
\[
\tilde{\epsilon}_l = \sum_{k=l}^\infty \epsilon_k {k\choose l} .
\]
In view of the assumption {\bf A1},
\begin{equation}
\label{tilde_epsilon}
\abs{\tilde{\epsilon}_l} \leq \delta \sum_{k=l}^\infty\rho^{-k}
{k\choose l}  = \delta \frac{\rho^{-l}}{(1-1/\rho)^{l+1}} =
\frac{\delta\rho}{(\rho -1 )}(\rho -1)^{-l} .
\end{equation}
Consequently,
\begin{equation}
\label{upper2}
\left|\bbV (\xi )
\right| \leq \frac{\delta\rho}{(\rho -1 )}
\sum_1^\infty \frac{\eta^{l-1}}{(\rho -1)^l} = \frac{\delta\rho}{(\rho -1)
(\rho - (1+\eta))} ,
\end{equation}
whenever $\xi\in \bbD_\eta (1)$.

\medskip
Pick $\eta = (\rho -1)/2$. Then the right hand sides of both
\eqref{upper1} and \eqref{upper2} are bounded above by
$2\delta (1+\rho)/(\rho -1 )^2$. Only at this stage we evoke
assumption {\bf A2}: under an appropriate choice of $\alpha$ the
latter expression is as small as desired, say less than $1/6$.
In view of \eqref{lower},\eqref{upper1} and \eqref{upper2} we,
therefore, conclude:
\begin{lem}
\label{lem:lower}
 Assume {\bf A1} and {\bf A2}, Then,
\begin{equation}
\label{eq:lower}
\min_{\xi \in \bbD_{(1+\rho)/2}(0)}
\left| \frac{\bbF (\xi )-1}{\xi -1}\right| \geq \frac13 .
\end{equation}
\end{lem}
Finally, by direct computation:
\begin{equation}
\label{Q_bound}
\left|\frac{\bbQ (\xi )-1}{\xi -1}\right| = \left|\frac1{1-p\xi}\right|
\geq \frac1{1+(1+\eta )p}\geq \frac12 .
\end{equation}

\subsection{Upper bound on
the numerators}
\label{sub:upper}
We continue to employ the notation of the preceeding subsection.
In particular,
\[
\bbV (\xi ) \df \frac{\bbE (\xi)}{1-\xi} \ \ {\rm and}\ \
\frac{\bbE (\xi)}{1-\xi} - \bbE^\prime (1) = \bbV (\xi ) - \bbV (1) .
\]
Since  by \eqref{eps_sum} , $\sum\epsilon_k =\bbE (1) = 0$, $\bbV$ is analytic
on $\bbD_{\rho}(0)$.
As in Subsection~\ref{sub:lower} pick $\eta\in (0,\rho -1)$ and
consider the following two cases:

\medskip
\noindent
\case{1}.  $\xi\in \bbD_{1+\eta}(0)\setminus\bbD_{\eta }(1)$. By
\eqref{upper1} and \eqref{upper2}
\[
\left| \frac{\bbV (\xi ) -\bbV (1)}{\xi -1} \right| \leq
\frac{\delta\rho}{\eta (\rho -1)
(\rho - (1+\eta))} +
\frac{\delta(1+\eta)}{\eta^2 (\rho -(1+\eta ))}\leq
\frac{2\delta \rho}{\eta^2 (\rho -(1+\eta ))}  .
\]
%

\medskip
\noindent
\case{2}. $\xi \in \bbD_{\eta }(1)$. Since $\bbV$ is analytic on
$\bbD_\eta (1)$ we , employing the notation of Subsection~\ref{sub:lower},
estimate:
\[
\left|\frac{\bbV (\xi )-\bbV (1)}{\xi -1}\right| \, =\,
\left|\sum_{l=2}^\infty \tilde{\epsilon}_l (\xi -1)^{l-2}\right|
\leq \frac{\delta\rho}{(\rho -1)^2 (\rho -(1+\eta ))} ,
\]
where we have performed a straightforward series summation
bounding $\abs{\tilde{\epsilon}_l}$ as in \eqref{tilde_epsilon}.

Picking $\eta =(\rho -1)/2$  and $R= (1+\rho)/2$ we
infer:
\begin{lem}
\label{lem:upper}
 Assume {\bf A1}, Then,
\begin{equation}
\label{eq:upper}
\max_{\xi\in \bbD_{R}}
\left|\frac{\bbE (\xi ) -(\xi - 1)\bbE^{\prime} (1)}{(\xi -1 )^2}
\right|
\, \leq \,
\frac{2\delta\rho}{(\rho -1)^3}
\end{equation}
\end{lem}
On the other hand,
\[
\frac1{\xi -1}\left( \frac{\bbQ (\xi )-1}{\xi - 1}- \bbQ^\prime (1)\right)
= \frac{p}{q(1-p\xi )}.
\]
Since by assumption {\bf A1}, $\max_{\abs{\xi }=R}\abs{\bbE (\xi )}\leq
\delta\rho/(\rho -1 )$, we arrive to the following bound for the
numerator of $\bbU_2$:
\begin{equation}
\label{U2_upper}
\max_{\abs{\xi }=R}
\left|\frac{\left(\bbQ (\xi) - 1 - \bbQ^\prime (1)(\xi -1)\right)\bbE (\xi)}
{(\xi -1)^2}\right| \leq \frac{p}{q(1-Rp)}\frac{\delta\rho}{(\rho - 1)}
\leq \frac{2\delta\rho}{p(\rho - 1)^3} .
\end{equation}

\subsection{The target bound on $\abs{r_n}= \abs{g_n -g}$}
\label{sub:target}
As before set  $R= (1+\rho )/2$.
By the estimates of Lemma~\ref{lem:lower} and
Lemma~\ref{lem:upper},
\[
\max_{\abs{\xi \leq R}} \abs{\bbU_1 (\xi )}
+\max_{\abs{\xi \leq R}} \abs{\bbU_2(\xi )}
 \leq
 \frac{6\delta\rho (2+p)}{p (\rho -1)^3} .
\]
Finally, $\bbF^{\prime}(1) =1/q +\sum k\epsilon_k$. By
assumption {\bf A1},
\[
\abs{\sum_k k\epsilon_k} \leq \delta \sum_1^\infty kp^k =
\frac{\delta p}{q^2}\leq \frac{\delta}{\rho -1}\frac1{q}.
\]
By the scaling relation between $\rho$ and $\delta$ (assumption {\bf A2})
the right hand side above is ${\small o}(q)$. In particular,
\begin{equation}
\label{gqcontrol}
 g  = q (1 + \small{o}(1)) ,
\end{equation}
as it now follows from \eqref{FN_limit}.

Substituting the above estimates into \eqref{exp_bound}:
\begin{thm}
\label{thm:convergence}
Assume {\bf A1} and {\bf A2}. Set $R= (1+\rho)/2$. Then,
\begin{equation}
\label{gnterm}
\begin{split}
\abs{r_n} &=
\abs{ g_n - g}\, \leq\, \frac{12 q\delta\rho (2+p)}{p (\rho -1)^3}
R^{-n} =
\frac{12 q\delta\rho (2+p)}{p (\rho -1)^3}
\left(\frac2{1 +\rho }\right)^n\\
&\df c_1 (q ,\rho )\delta\left(\frac2{1 +\rho }\right)^n .
\end{split}
\end{equation}
In particular for every $\nu <(\rho -1)/4$,
\begin{equation}
\label{nu_sum}
\begin{split}
\sum_n \abs{r_n} (1+\nu )^n &=
\sum_n \abs{g -g_n } (1+\nu )^n
\leq \frac{12\rho (2+p)(1+\rho)}{p}
\frac{q\delta}{(\rho -1- 2\nu )(\rho -1)^3}\\
&\leq \frac{48\rho (2+p )(1+p)}{p (\rho - 1)^4} \delta g
  \,  \df\,  c_2 (q,\rho ) \delta g .
\end{split}
\end{equation}
\end{thm}

\subsection{Uniform rod activities}
\label{sub:uniform}
Let the rod activities be given by
\begin{equation}
\label{eq:unif_rod}
f_k = \left\{
\begin{split}
&\bar{q}\, ;\, k=1,\dots ,N \\
&0\, ;\, {\rm otherwise}
\end{split}
\right.
\end{equation}
Above $\bar{q} = (\sum_1^N p^{k-1})^{-1} = (1-p)/(1-p^N )$. Thus,
$\bar{q} - q = qp^N /(1-p^N )$. In other words the sequence of weights
$ \{f_k\}$ in \eqref{eq:unif_rod} corresponds, in the notation
of \eqref{f_type}, to
\begin{equation}
\label{eps_uniform}
\epsilon_k = \left\{
\begin{split}
&\frac{qp^N}{1-p^N} p^{k-1}\, ;\, k\leq N\\
&-qp^{k-1}\,\ \ \ \ \  ;\, k>N
\end{split}
\right.
\end{equation}
Without loss of generality we may assume that $q<1/2$. Then for each
$N$ fixed the weights $\{ \epsilon_k \}$ satisfy assumption {\bf A1}
with , for example,
\begin{equation}
\label{choices}
 \rho = \left( 1+ \frac1{p}\right)/2 = 1+ \frac{q}{2(1-q)}\quad
 \text{and} \quad \delta = \delta_N (q) = \left(
 1-\frac{q}{2}\right)^N .
 \end{equation}
 Of course, assumption {\bf A2} will be
 also satisfied for such choice of $\rho$ and $\delta$ as soon
as
\[
\delta_N (q) =(1- \frac{q}{2})^N \leq \alpha (\rho -1)^2 = \alpha
\left( \frac{q}{2(1-q)} \right)^2 .
\]
For the value of $\rho$ related to $q$ as in \eqref{choices} set
\[
\bar{c}_1 (q) =c_1 (q, \rho )\quad \text{and}\quad \bar{c}_2 =c_2
(q,\rho ) ,
\]
where $c_1$ and $c_2$ are the universal constant which appear on
the right hand sides of \eqref{gnterm} and \eqref{nu_sum}. Let us
reformulate the claim of Theorem~\ref{thm:convergence} as applied
to the case of uniform rod activities (with the scaling choice
\eqref{choices} in mind):

\begin{lem}
\label{lem:milos_bound}
Let $q<1/2$ be fixed.
Then there exists $N_0 = N_0 (q)$
 such that for every $N\geq N_0$,
 the uniform rod weights
$\{f_k\}$ in \eqref{eq:unif_rod} generate
the renewal sequence $\{ g_n\}$ which satisfies:
\begin{equation}
\label{Brn} \abs{r_n} =\abs{g_n -g} \leq \bar{c}_1 (q)\delta_N (q)
\left( 1+ \frac{q}{4(1-q)}\right)^{-n} ,
\end{equation}
where $g$ was defined in \eqref{FN_limit}. Moreover, for $\nu
=(\rho -1)/4 = q/8(1-q)$,
\begin{equation}
\label{Bsum} \sum_{n} \abs{r_n}
(1+\nu )^n \, \leq \,
\bar{c}_2 (q ) \delta_N (q ) g .
\end{equation}
\end{lem}

\section{Perturbation theory}
\label{ssec_PS}
In this section $V$ is the lattice torus of a fixed (large) linear size
$L$; $V={\mathbb Z}^2/\text{mod}(L)$. Notice, however, that all the 
estimates below do not depend on $L$.

\subsection{Super-contours}
We proceed similarly as in the proof of Theorem~\ref{thm_exact}. We first
 split vacancies into two families, and partition the box 
$V$ 
into the two disjoint sub-boxes $V_{\mathrm{h}}$  and $V_{\mathrm{v}}$
containing the horizontal, resp. vertical, sites.  Associated to this
partition, there is a family of one-dimensional   
boxes $\underline{\Delta} = (\Delta_i)$, each of
 which is either a horizontal ``segment'' in $V_{\mathrm{h}}$, or a vertical
 ``segment'' in $V_{\mathrm{v}}$. The weight of 
the partition can then be expressed as a product
 over all $\Delta \in \underline{\Delta}$ of the
 corresponding one-dimensional partition
 functions $g_{|\Delta|} = Z^{\mathrm{1D}}_{|\Delta|}$.
 Contrarily to what happens in the case considered in
 Theorem~\ref{thm_exact}, these one-dimensional
 partition function \emph{do} generally depend on the
 length of the corresponding box $\Delta$. However, 
as we have seen in Section~\ref{asymptotics1DPF},
 these partition functions approach their limiting
 value $g$ rather quickly, provided we choose $N$ 
large enough. It is therefore convenient to expand them
 around this limiting value:
$$
\prod_{\Delta\in\underline{\Delta}} g_{|\Delta|} =
 \prod_{\Delta\in\underline{\Delta}} g 
\left( 1 + \frac{g_{|\Delta|}-g}{g} \right) \,.
$$
We want to use this expansion in order to obtain a perturbation of the
 pure Ising model which appeared in the case of Theorem~\ref{thm_exact}.
 Let us denote by $\underline{\gamma}=(\gamma_i)$ the family of Ising 
contours appearing when interpreting $V_{\mathrm{h}}$, resp. $V_{\mathrm{v}}$,
 as the region occupied by $+$, resp. $-$, spins. We can then associate 
to each of these contours a weight $w(\gamma) = e^{-2\beta|\gamma|}$, 
where we have set $e^{-2\beta} = \sqrt{g}$; this allows us to write simply
$$
\prod_{\Delta\in\underline{\Delta}} g = \prod_{\gamma\in\underline\gamma}
 e^{-2\beta|\gamma|}\,.
$$
We would like to encode all the information from the partition 
into these contours; in order to do this, we suppose that 
these contours come with a ``color'', i.e. each contour 
$\gamma$ carries the information on which of the two sets
 $V_{\mathrm{h}}$ or $V_{\mathrm{v}}$ belong to which side
 of the contour. Of course, there is then a compatibility 
condition on these contours (in addition to their being 
disjoint): the colors must match.

We also introduce the set of excited intervals $\underline I = (I_i) 
\subset \underline\Delta$, and associate to such objects the weight
 $w(I) = (g_{|I|}-g)/g$. Using this we can write
$$
\prod_{\Delta\in\underline{\Delta}} 
\left( 1 + \frac{g_{|\Delta|}-g}{g} \right) =
 \sum_{\underline I \subset \underline\Delta} \prod_{I\in\underline I} w(I)\,.
$$
Of course, since our colored contours $\underline\gamma$ contain 
all the information on the partition, the family $\underline\Delta$ 
is actually completely determined by the contours.

We now introduce our basic notion of super-contours, which are maximal 
connected components of (colored) contours and excited intervals
 (saying that an interval is connected to a contour if at least 
one of the extremities of the interval belongs to the contour).
 We denote the family of super-contours by $\underline\Gamma=(\Gamma_i)$. 
The weight $w(\Gamma)$ of a super-contour $\Gamma$ is then naturally 
given by the product of the weights of the contours and excited 
intervals it encompasses. We therefore finally obtain the following
 expression for the weight of the total partition function of our model:
$$
Z_{q,N,L} = \sum_{\underline\Gamma} \prod_{\Gamma\in\underline\Gamma} 
w(\Gamma)\,,
$$
where the sum is taken over all compatible families of  
super-contours, i.e. those resulting from a partition
 $V=V_{\mathrm{h}}\vee V_{\mathrm{v}}$ in the way just described.
\begin{figure}[t!]
   \centering
   \scalebox{0.3}{\input{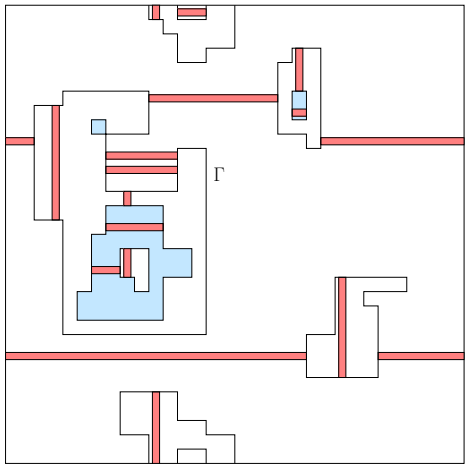}}
   \caption{A configuration of super-contours on the torus $V$. 
There are four super-contours. Notice that $\Gamma$ if of
 horizontal type, but also possess some 
interior components of horizontal type (shaded in the picture). There are two super-contours winding around the torus; the contribution of configurations containing such super-contours being negligible when the box is large, as shown in Subsection~\ref{ssec_long}, they can actually be neglected.}
   \label{fig:supercontours}
\end{figure}

At this stage, it is not possible to apply a simple 
Peierls argument in order to control our model. 
Indeed, our super-contours are colored, and even 
though there is a symmetry in our model (under a 
simultaneous rotation by $\pi /2$  and exchange 
of horizontal and vertical sites.
On the other hand, there is also a fundamental asymmetry:
The shape of a region generally strongly favours one of the two species.
 This forces us to use the general strategy of the Pirogov-Sinai theory,
 which turns out to be quite simple in our case, due to the fact that,
 because of the above-mentioned symmetry, the free energies of the two
 phases are necessarily equal, and thus we are not required to add a 
suitable external field to reach phase coexistence.

The basic idea of the Pirogov-Sinai theory is to expand the 
partition function only over external contours, and 
introduce new weights, which reduce the compatibility condition to
 something of purely geometrical nature. 

However, since we work on the lattice torus 
$V ={\mathbb Z}^2/\text{mod} (L)$, the notion of exteriour of 
a contour  is ambiguous. One way to mend the situation would be to fix a 
distinguished site, say 0, and to declare it to be ``a point at infinity''.
 On the other hand all our computations below are based on relatively crude 
combinatorial estimates which take into account local graph geometry 
of ${\mathbb Z}^2$, but not the global topological structure of $V$. 
Consequently, we shall from the start ignore (necessarily long) winding
super-contours and then simply notice  that should we use the
``point at infinity'' definition of exteriour, the analog of \eqref{G_length}
 below would anyway render long winding contours improbable.

For any non-winding super-contour $\Gamma$ the exteriour of 
$\Gamma$ is defined in a straightforward fashion and, accordingly, 
the type of a non-winding super-contour $\Gamma$ will be declared to be 
horizontal or vertical if such is the colour of
its exteriour.
 
Thus, let ${\mathcal S}_L^h$ (respectively ${\mathcal S}_L^v$ )
  be the set of all non-winding horizontal (respectively vertical) type
  super-contours
on $V= {\mathbb Z}^2/\text{mod}(L)$. Of course, ${\mathcal S}_L^h$ and 
${\mathcal S}_L^h$ are related by the $\pi/2$ rotation symmetry: If 
$\Gamma\in {\mathcal S}_L^h$, then $\theta_{\pi/2}
\Gamma\in {\mathcal S}_L^v$. The interiour
 $\mathrm{int}(\Gamma)$  of coloured 
super-contours $\Gamma\in {\mathcal S}_L^h\cup  {\mathcal S}_L^v$
 is also coloured and in the sequel we shall write
\[
\mathrm{int}(\Gamma) = \mathrm{int}_{\mathrm{h}}(\Gamma) 
\cup \mathrm{int}_{\mathrm{v}}(\Gamma) 
\]
for the horizontal and vertical parts of $\mathrm{int}(\Gamma)$. 
%
By the $\pi /2$-rotation symmetry of $V$, restricting to 
external contours of horizontal type (``h-type''), i.e. 
with their exterior colored as horizontal, yields 
exactly one-half of the full partition function: 
Using the notation $\underline{\mathcal S}_L^{h, \rm{ext}}$ for set of 
all collections of compatible {\em external} contours from ${\mathcal S}_L^{h}$
and, accordingly, 
$\underline{\mathcal S}_L^{h}$ for 
for set of 
all collections of compatible contours from ${\mathcal S}_L^{h}$, 
   we can then write
\begin{align*}
Z_{q,N,L}
&= 2\sum_{\underline\Gamma \in  \underline{\mathcal S}_L^{h, \rm{ext}}} 
\prod_{\Gamma\in\underline\Gamma} w(\Gamma) \, 
Z^{\mathrm{h}}_{\mathrm{int}_{\mathrm{h}}(\Gamma)} 
Z^{\mathrm{v}}_{\mathrm{int}_{\mathrm{v}}(\Gamma)}\\
&= 2\sum_{\underline\Gamma \in 
\underline{\mathcal S}_L^{h, \rm{ext} } }
\prod_{\Gamma\in\underline\Gamma} \widetilde{w}(\Gamma) \, 
Z^{\mathrm{h}}_{\mathrm{int}(\Gamma)}\\
&= 2\sum_{\underline\Gamma \in\underline{{\mathcal S}}_L^h } 
\prod_{\Gamma\in\underline\Gamma} \widetilde{w}(\Gamma)\,,
\end{align*}
where the new weights are given by
\begin{equation}
\label{eq:newWeights}
\widetilde{w}(\Gamma) = w(\Gamma) \,
\frac{Z^{\mathrm{v}}_{\mathrm{int}_{\mathrm{v}}(\Gamma)}}
{Z^{\mathrm{h}}_{\mathrm{int}_{\mathrm{v}}(\Gamma)}} .
\end{equation}
Notice that  in the last expression the sum is over 
{\em all} compatible families of  h-type super-contours; in particular, the 
compatibility condition is now purely geometrical.

\subsection{Cluster expansion}
The next step is to show that the 
new weights are still under control.
We first need a bit of terminology. Let us 
denote by $\Gamma\not\sim\Gamma'$ the 
relation ``$\Gamma$ is incompatible with $\Gamma'$''. A 
cluster is a family $\mathcal C$ of super-contours 
which cannot be split into two disjoint families 
$\mathcal{C}_1$ and $\mathcal{C}_2$ such that 
all pairs $\Gamma\in\mathcal{C}_1$ and $\Gamma'\in\mathcal{C}_2$ are
compatible. We also write $\mathcal C \not\sim \Gamma$ if 
there exists $\Gamma'\in\mathcal{C}$ such that 
$\Gamma'\not\sim\Gamma$. Finally we write $|\Gamma|$ 
for the total length of all the contours and intervals forming $\Gamma$.
We want to be able to use the following classical 
sufficient condition for the convergence of the cluster expansion~\cite{KP86}:
\begin{lem}
\label{lem:KP}
Suppose that, for some small $a>0$,
\begin{equation}
\label{sum_condition}
\sum_{\Gamma'\,:\, \Gamma' \not\sim \Gamma} e^{2a\,|\Gamma'|}
 \abs{\widetilde{w}(\Gamma')} \leq a\, |\Gamma|\,,
\end{equation}
for each $\Gamma$. Then $Z_{q,N,L} \neq 0$ and there 
exists a unique function $\Phi^\mathrm{T}$ on the set of clusters such that
$$
\log Z_{q,N,L} = \sum_{\mathcal{C} \subset V} \Phi^\mathrm{T}(\mathcal{C})\,.
$$
Moreover,
\begin{equation}
\label{cluster_a}
\sum_{\mathcal{C} \not\sim \Gamma}
 |\Phi^\mathrm{T}(\mathcal{C})| e^{a\, \abs{\mathcal{C}}} \leq a\, |\Gamma|\,.
\end{equation}
\end{lem}

We claim that the weights $\widetilde{w}$ indeed satisfy 
\eqref{sum_condition} once $q$ is chosen to be small enough and then 
$N$ is chosen sufficiently large. The argument comprises two steps: First we 
shall check that \eqref{sum_condition} holds for the weights $w$. Next we
shall argue that the conclusion of Lemma~\ref{lem:KP} for the weights
 $w$ actually implies the validity of \eqref{sum_condition} for
the target weights $\widetilde{w}$ for a possibly smaller value of
 $q$ and larger values of $N$. 

\begin{lem}
\label{Gamma_sum}
There exist $a>0$, $q>0 $ and $N_0 = N_0 (q)$ such that
\begin{equation}
\label{Gsum_condition}
\sum_{\Gamma'\,:\, \Gamma' \not\sim \Gamma} e^{2a\, |\Gamma'|}
 \abs{{w}(\Gamma')} \leq a\, |\Gamma|\,,
\end{equation}
for every $N\geq N_0$.
\end{lem}

\begin{lem}
\label{TGamma_sum}
There exist ${a}>0$, $\widetilde{q}>0 $ and $\widetilde{N}_0 
= \widetilde{N}_0 (\widetilde{q} )$ such that
\eqref{sum_condition} holds for $\widetilde{q}$ and 
for every $N\geq \widetilde{N}_0$.
\end{lem}

\subsection{ Proof of Lemma~\ref{Gamma_sum}} 
A  convenient way to over-count
$$
\sum_{\Gamma'\,:\, \Gamma' \not\sim \Gamma} e^{2a\,|\Gamma'|} \abs{w(\Gamma')}
$$
is as follows: Pick 
\begin{equation}
\label{achoice}
2a = \frac{\rho -1}{16} = \frac{q}{32 (1-q)} .
\end{equation}
 Any excited interval $I$ of a super-contour $\Gamma$ connects
two
dual bonds $b$ and $b^\prime$ which belong to Ising contours $\gamma$ and,
accordingly,
$\gamma^\prime$. There are two cases:

1) If $\gamma =\gamma^\prime$, we erase $I$ and upgrade  the weights of $b$ and
   $b^\prime$ from $\sqrt{g}$ to
\[
 \sqrt{g}e^{2a}  + \sum_{k\geq 1} e^{2a\, k} r_k .
\]

By \eqref{Bsum} and \eqref{choices} the latter expression is bounded 
above in absolute value by $2e^{-2\beta }$. 

2) If $\gamma \neq \gamma^\prime$, then we erase $I$ and instead add 
 two red links which
connect between the endpoints of $b$ and $b^\prime$. Precisely, if $b= (u,v)$
and
$b^\prime = (t,s )$ where the both pairs $\{ u,v\}$ and $\{t, s\}$ of dual
vertices are
recorded in the lexicographical order, then we add red links $(u,t )$ and
$(v,s)$. In this way both red links lie on the dual lattice and have the
same length $k=1,2,\dots $. We associate the weight $r_k $ to each of
those links.

Clearly after the above procedure is applied to all the excited intervals of
$\Gamma$ we end up with a connected edge 
self-avoiding polygon $\widehat{\Gamma }$
which entirely lies on the dual lattice. In order to control the original
weights
${w} (\Gamma )$ we over-count via ignoring the
geometric constraints: from each vertex of the dual lattice one is permitted
to
grow up bonds in all 4 possible directions: either usual Ising bonds with
weights  $e^{-2\beta}$ or ``red'' bonds of lengths $k=1,2,\dots$ with
the weights $e^{2ak} r_k$ respectively. Notice that any modified graph
$\widehat{\Gamma}$ contains at least 4 Ising bonds. Consequently, 
\begin{equation}
\label{four_sum}
\begin{split}
\sum_{0\in{\Gamma}}& e^{2a\, |\Gamma |}\abs{w (\Gamma )} \\
& \leq
\sum_{n=4}^{\infty} 4^n \left( 2e^{-2\beta} +\sum_1^\infty
  ke^{2a\, k} r_k
 \right)^n ,
\end{split}
\end{equation}   
where the above sum is over the total number of bonds (and links ) of
$\widehat{\Gamma }$.  By \eqref{Brn} and \eqref{Bsum} and in view of the 
possibility to control the smallness of $\delta_N$ via \eqref{choices}, 
we, given small $q$ and $a$ as in \eqref{achoice}, can always choose
 a large enough value of $N_0$, such that,
\[
\left| 2e^{-2\beta  } +\sum_1^\infty
  ke^{2a\,  k} r_k\right| \leq  3 e^{-2\beta  } .
\]
As a result the right hand side of \eqref{four_sum} is bounded above by
\[
\frac{(12 e^{-2\beta  } )^4}{1 - 12 e^{-2\beta  }} .
\]
Recalling the notation $e^{-2\beta  } = \sqrt{g}$ we, in view of 
\eqref{gqcontrol},  infer that the latter expression 
  is much less than the value of $a$ in \eqref{achoice} once 
$q$ happens to be  sufficiently small.


\vskip 0.2cm

In the sequel we shall assume that $q\leq q_0$ and $N\geq N_0 (q)$ 
are such that we actually have a strengthened version 
of \eqref{new_condition}: Set $w_0 (\Gamma ) =\abs{w (\Gamma )}$.
 Then,  
\begin{equation}
\label{Gamma_length}
\sum_{{\Gamma}\ni 0}\abs{\Gamma} e^{2a \abs{\Gamma}} w_0 (\Gamma )
\leq a.
\end{equation}
Indeed, \eqref{Gamma_length} follows by a straightforward adjustment of the 
arguments employed for the proof of Lemma~\ref{Gamma_sum}.

\subsection{ Proof of Lemma~\ref{TGamma_sum}}
Let $q$ and $N_0$ are fixed as in the proof of Lemma~\ref{Gamma_sum},
 and let $\{ w_0 (\Gamma )\df \abs{w (\Gamma )}\}$ are the 
the absolute values of the weights of super-contours
 evaluated at such values of $q$ and $N_0 (q )$. It is enough to check
that there is $\widetilde{q}\leq q$ and $\widetilde{N}_0\geq N_0$, such that
for every $N\geq \widetilde{N}_0$, the $(\widetilde{q} ,N)$ super-contour
weights $\{ (w (\Gamma )\}$ satisfy:
\begin{equation}
\label{new_condition}
 \abs{w (\Gamma )} e^{\sum_{\gamma\in\Gamma }\abs{\Gamma}} \leq 
w_0 (\Gamma )\quad \text{and}\quad  
\frac{Z^{\mathrm{v}}_{V}}{Z^{\mathrm{h}}_{V}} \leq e^{|\partial V|}\,
\end{equation}
for every super-contour $\Gamma$ and for each finite subset 
$V\subset {\mathbb Z}^2$. 

Of course, only the second inequality in 
\eqref{new_condition} deserves  to be checked.  
 This is 
done by induction on the volume. 
 Obviously, if the volume $|V|=1$, then we have
$
Z^{\mathrm{v}}_{V} / Z^{\mathrm{h}}_{V} = 1\,.
$
Suppose now that indeed 
$$
\frac{Z^{\mathrm{v}}_{V}}{Z^{\mathrm{h}}_{V}} \leq e^{|\partial V|}\,,
$$
for all $\abs{V} < K$.
We want to prove that this also holds when $|V|=K$. In order 
to see that, observe that all the super-contours appearing in these two
partition functions have interiors of volume at most $K-1$. Introducing 
the sets $\mathcal{S}_{K-1}^{\mathrm h}$ and 
$\mathcal{S}_{K-1}^{\mathrm v}$ of all clusters 
made up of h-type, resp. v-type, super-contours having (total) 
interior of volume at most $K-1$, and using the symmetry 
present in the model, we can write
$$
\frac{Z^{\mathrm{v}}_{V}}{Z^{\mathrm{h}}_{V}} = 
\frac{Z^{\mathrm{v}}_{V}}{\exp\Bigl(\sum_{x\in V}\sum_{
\substack{\mathcal{C}\in\mathcal{S}_{K-1}^{\mathrm v}\\
 \mathcal{C}\ni x}}|\mathcal C\cap V|^{-1}\,
\Phi^{\mathrm{T}}(\mathcal{C})\Bigr)} \frac{\exp\Bigl(\sum_{x\in V}
\sum_{\substack{\mathcal{C}\in\mathcal{S}_{K-1}^{\mathrm h}\\ 
\mathcal{C}\ni x}}|\mathcal C\cap V|^{-1}\,
\Phi^{\mathrm{T}}(\mathcal{C})\Bigr)}{Z^{\mathrm{h}}_{V}}\,.
$$
Notice now that all the contours $\Gamma$  appearing in the
 above partition functions have weights $\widetilde{w} (\Gamma )$
 which, by the induction
assumption and by the first of the inequalities in 
\eqref{new_condition}, 
 satisfy:
$$
\abs{\widetilde{w}(\Gamma)}
 \leq \abs{w(\Gamma)}\, e^{\sum_{\gamma\in\Gamma}|\gamma|}\leq
w_0 (\Gamma ) .
$$
Therefore 
we can apply Lemma~\ref{lem:KP}. Expanding the two partition 
functions and cancelling the terms involving clusters 
entirely contained inside $V$, we obtain the desired result, since
by \eqref{cluster_a}
$$
\frac{Z^{\mathrm{v}}_{V}}{Z^{\mathrm{h}}_{V}}  
\leq e^{2a \, |\partial V|}\,,
$$
and $2a <1$ once, according to \eqref{achoice}, $q$ is not very close
to $1$. 
 

\section{Proofs of the main results}
In this section we complete the proofs of Theorem~\ref{thm_Main} 
and Theorem~\ref{thm_Domains}. As in the proof of Lemma~\ref{TGamma_sum}
we proceed to work within the range of parameters $({q} ,N)$
 which satisfy \eqref{new_condition}

\subsection{Contribution of long super-contours}
\label{ssec_long}
As before let $V$ be a lattice torus of linear size $L$. Given 
a supercontour $\Gamma\in{\cal S}_L^h$ define 
\begin{equation*}
\begin{split}
G_{q,N,L} (\Gamma )\,  &= \,
\frac{1}{Z_{q,N,L }}\sum_{\underline{\Gamma}\ni\Gamma}
\widetilde{w}( \underline{\Gamma} )\\
& = \, 
\frac12 \widetilde{w}(\Gamma ) exp\left(\, - 
\sum_{\mathcal{C} \not\sim \Gamma}
 \Phi^\mathrm{T}(\mathcal{C})\right) .
\end{split}
\end{equation*}
By \eqref{new_condition} and \eqref{cluster_a},
\[
\left| G_{q,N,L} (\Gamma ) \right| \leq w_0 (\Gamma ) e^{a\abs{\Gamma}} .
\]
Furthermore, by \eqref{sum_condition}, there exist  
constants $c_1$ and $c_2$ such that 
\begin{equation}
\label{G_length}
\begin{split}
\sum_{\abs{\Gamma}\geq k}
w_0 (\Gamma )e^{a \abs{\Gamma}}\, &\leq\, c_1 L^2 e^{-ak}
\sum_{\Gamma\ni 0} w_0 (\Gamma )e^{2a \abs{\Gamma}}\\
&\leq\, c_2 L^2 a e^{-ak} .
\end{split}
\end{equation}
As a result, there exists $c_3 <\infty$, such that the contribution 
of 
super-contours $\Gamma\in {\cal S}_{L}^h$ with 
$\abs{\Gamma} >c_3 \log L$  to the  partition function
$Z_{q,N,L }$  is, uniformly in $L$, negligible. The same argument
applies, of course, in the case of vertical super-contours
 $\Gamma\in {\cal S}_L^v$.

\subsection{Proof of Theorem~\ref{thm_Main}}
The first statement of Theorem~\ref{thm_Main} follows immediately from
results of Gruber and Kunz~\cite{GrKu71}.

Assume now that the 
parameters $({q} ,N)$
  satisfy \eqref{new_condition}. By \eqref{Gamma_length} we may
exclude long winding contours. Thus, 
the only thing remaining to be done in order to complete the proof
of Theorem~\ref{thm_Main} is to estimate the probability that 
a given site, say $0$, belongs to the interior of some short 
non-winding contour. In view of Lemma~\ref{lem:KP} the probability that
$0$ belongs to the interior of such a super-contour can then be written as
\begin{equation}
\label{prob_zero}
\sum_{\Gamma \circlearrowright 0} w (\Gamma)
\frac{Z^{\mathrm{v}}_{\mathrm{int}_{\mathrm{v}}(\Gamma)}
Z^{\mathrm{h}}_{\mathrm{int}_{\mathrm{h}}(\Gamma)}}{Z^{\mathrm{v}}_{q,N,V}}\, 
=\, \sum_{\Gamma \circlearrowright 0} \widetilde{w} (\Gamma)
\exp\left(-\sum_{C\not\sim\Gamma} \Phi^{\mathrm{T}}(\mathcal{C}) \right) .
\end{equation}
where $\Gamma \circlearrowright 0$ means that $0$ is in the interior of the
super-contour $\Gamma$. By \eqref{new_condition} the latter expression 
is bounded above by  
\[
\sum_{\Gamma \circlearrowright 0} w_0 (\Gamma ) e^{a\abs{\Gamma}} 
\leq \sum_{\Gamma \ni 0}
\abs{\Gamma} w_0 (\Gamma ) e^{a\abs{\Gamma}} ,
\]
 The claim of the Theorem follows now from 
\eqref{Gamma_length}.

\subsection{Infinite volume states} Let ${\cal A}_\infty$ be the set of all
such coverings $\widetilde{\omega }$ 
 of ${\mathbb Z}^2$ by horizontal and vertical rods (we colour 
monomers as well), which contain only finite contours. Of course, for every 
$\widetilde{\omega}$ the notion of the exteriour colour 
$\chi (\widetilde{\omega}) = h\ \text{or}\ v$ is well defined. By a
 straightforward application of Lemma~\ref{lem:KP}:
\begin{thm}
\label{thm:unique}
There exists $q_0 > 0$ such that for every $q\leq q_0$ one can find 
$N_0 = N_0 (q )$ which enjoys the following property: For every 
$N\geq N_0$ there exists a unique infinite volume Gibbs state $\mu_{q,N}^h$
(respectively $\mu_{q,N}^v$) such that
\[
\mu_{q,N}^h\left( {\cal A}_\infty ;\chi (\widetilde{\omega })= h\right)\ =\ 1
\]
(respectively 
$\mu_{q,N}^v\left( {\cal A}_\infty ;\chi (\widetilde{\omega })= v\right) 
= 1$) .
Furthermore, let $\underline{\gamma}$ be the ({\em random}) set of of all the 
exteriour contours of $\widetilde{\omega}$ and, given a finite domain 
$\Lambda\subset {\mathbb Z}^2$, let $\underline{\gamma}_\Lambda = 
(\gamma_1 ,\dots ,\gamma_n)$ be a {\em fixed} 
compatible set of exteriour contours, 
such that each $\gamma_k$ intersects $\Lambda$, $\Lambda\cap\gamma_k\neq
\emptyset$. Then, 
\begin{equation}
\label{muh_formula}
\mu_{q,N}^h \left(\underline{\gamma}_\Lambda \subset\underline{\gamma}\right)
\, =\, 
\sum_{\underline{\Gamma}_\Lambda\sim\underline{\gamma}}
\widetilde{w}({\underline{\Gamma}_\Lambda}) \exp\left( - 
\sum_{{\mathcal C}\not\sim \underline{\Gamma}_\Lambda} \Phi^{\mathrm T}
({\mathcal C})\right) ,
\end{equation}
where the above sum is over all compatible collections 
$\underline{\Gamma}_\Lambda = (\Gamma_1 ,\dots ,\Gamma_m)$ of 
super-contours satisfying:
\[
\forall\, l=1,\dots ,m\,\exists\, k\ \text{such that}\ \gamma_k\in\Gamma_l
\quad\text{and}\quad \cup\gamma_k\subseteq\cup\Gamma_l .
\]
\end{thm}
Formulas \eqref{muh_formula} and \eqref{cluster_a} readily imply that 
 $\mu_{q,N}^h$ has an exponential clustering property: Given two
 disjoint  
boxes $\Lambda_1$ and $\Lambda_2$ and two fixed compatible collections
 $\underline{\gamma}_{\Lambda_1}$ and $\underline{\gamma}_{\Lambda_2}$ with 
 $\underline{\gamma}_{\Lambda_k}\subseteq \Lambda_k ;\ k=1,2$, the 
following 
bound holds:
\begin{equation}
\label{clustering}
\left| \log\frac{
\mu_{q,N}^h\left(\gamma_{\Lambda_1}\subset\underline \gamma\ ;
 \gamma_{\Lambda_2}\subset\underline \gamma\right)}
{\mu_{q,N}^h\left(\gamma_{\Lambda_1}\subset\underline \gamma\right)
\mu_{q,N}^h\left(\gamma_{\Lambda_2}\subset\underline \gamma\right)
}
\right|\ \leq\ c_1\abs{\Lambda_1}\abs{\Lambda_2} e^{-ac_2 {\mathrm d}
(\Lambda_1 ,\Lambda_2 )},
\end{equation}
where ${\mathrm d}
(\Lambda_1 ,\Lambda_2 )$ is a distance (say $l_1$) and $c_1$ and $c_2$
are two positive constants which depend only on $q$ and $N$. 

\subsection{Boundary surface tension}
Consider vertical and horizontal intervals
\[
J_k^v = (1/2, 1/2) + \{ (0,0), (0,1 ),\dots ,(0,k-1)\}\quad\text{and}
\quad J_k^h = (1/2, 1/2) + \{ (0,0),\dots (k-1 ,0)\} .
\]
By construction, both $J_k^v$ and $J_k^v$ are linear segments on the 
dual lattice $(1/2, 1/2) +{\mathbb Z}^2$. Given a rod
 $I = (u_1,\dots ,u_n)\subset {\mathbb Z}^2$ let us say that 
$I$ intersects  $J_k^v$; $I\cap J_k^v\neq \emptyset$ if
\[
J_k^v\cap \text{int}\left(\cup_{k=1}^n B_1 (u_k )\right)\neq \emptyset ,
\]
where $B_1 (u) = u+ [-1/2 ,1/2]\times[-1/2 ,1/2]\subset {\mathbb R}^2$ and 
for a bounded set $A\subset {\mathbb R}^2$ the symbol 
$\text{int}(A)$ stands for  its ${\mathbb R}^2$-interiour. In a similar 
fashion we define $I\cap J_k^h\neq \emptyset$.
  Notice that monomers cannot intersect $J_k^v$ or 
$J_k^h$. Also, with such a definition, 
$J_k^v$ cannot be intersected by a vertical rod and, accordingly, 
$J_k^h$ cannot be intersected by a horizontal one. 

Given  a ${\mathbb Z}^2$ tiling $\widetilde{\omega}\in {\cal A}_\infty$
let us say that the event $\{ J_k^v\cap\widetilde{\omega }=\emptyset\} $ (
respectively $\{ J_k^h\cap\widetilde{\omega }=\emptyset\} $) occurs if 
$J_k^v$ (respectively $J_k^h$) does not intersect any of the rods 
of $\widetilde{\omega }$. 

\vskip 0.2cm
\noindent
We define two types of boundary surface tensions:
\begin{equation}
\label{tau_qN}
\tau_{q,N}\, =\, -\lim_{k\to\infty}\frac1k\log \mu_{q,N}^h\left(
 J_k^v\cap\widetilde{\omega }=\emptyset\right)\, =\, 
-\lim_{k\to\infty}\frac1k\log \mu_{q,N}^v\left(
 J_k^h\cap\widetilde{\omega }=\emptyset\right) ,
\end{equation}
and 
\begin{equation}
\label{xi_qN}
\xi_{q,N}\, =\, -\lim_{k\to\infty}\frac1k\log \mu_{q,N}^h\left(
 J_k^h\cap\widetilde{\omega }=\emptyset\right)\, =\, 
-\lim_{k\to\infty}\frac1k\log \mu_{q,N}^v\left(
 J_k^v\cap\widetilde{\omega }=\emptyset\right) ,
\end{equation}
In both cases the fact that the corresponding quantities are well 
defined follows from standard sub-additivity arguments based
 on the exponential clustering property \eqref{muh_formula} and on the 
$\pi/2$-rotational symmetry between the vertical and horizontal states.

\begin{lem}
\label{hv_gap} For any $q$ sufficiently small there exists $N_0 =N_0 (q)$
such that for every $N\geq N_0$,
\begin{equation}
\label{eq:hv_gap}
\tau_{q,N}\,  >\, \xi_{q,N} .
\end{equation}
\end{lem}
\noindent{\em Proof.} Since we do not try to prove the lemma in the 
whole range of entropy driven symmetry breaking, the poof boils down to
 a crude
perturbative argument. We start with a lower bound on $\tau_{q,N}$: Fix 
$\gamma_1,\dots ,\gamma_n$ to be the set of all exteriour contours of 
$\widetilde{\omega}$ which intersect $J_k^h$. By \eqref{gqcontrol} and 
\eqref{gnterm}, 
\begin{equation} 
\label{Bound1}
\mu_{q,N}^v\left( \widetilde{\omega}\cap J_k^h=\emptyset\, \Big|\,
\gamma_1 ,\dots ,\gamma_n\right)\, \leq\, (2q)^{\abs{ J_k^h\setminus \cup 
\gamma_l}}.
\end{equation}
It remains, therefore, to derive an upper bound on 
\[
\mu_{q,N}^v\left(\abs{ J_k^h\setminus\cup_l
\gamma_l} \leq k/2\right) .
\]
By a straightforward modification of the over-counting argument employed
in the proof of Lemma~\ref{Gamma_sum} we infer from \eqref{muh_formula}
that for a given collection $\gamma_1 ,\dots, \gamma_n$ of 
exteriour contours,
\[
\mu_{q,N}^v\left( \gamma_1 ,\dots ,\gamma_n \right)\, \leq\, 
\exp\left( -2\beta \sum_1^n \text{diam} (\gamma_l )\right) ,
\]
where, as before, $e^{-2\beta } = \sqrt{g}$ and $g$ is related to $q$
 via \eqref{gqcontrol}. Elementary combinatorics leads then to 
the following conclusion: If $q$ is sufficiently small and 
$N\geq N_0 (q)$, then 
\begin{equation}
\label{Bound2}
\mu_{q,N}^v\left(\abs{ J_k^h\setminus\cup 
\gamma_l} \leq k/2\right) \, \leq\, \exp\left( -\frac{\beta k}{2}\right) .
\end{equation}
Combining \eqref{Bound1} and \eqref{Bound2} we arrive to the following
lower bound on $\tau_{q, N}$:
\begin{equation}
\label{tau_bound}
\tau_{q, N} \,  \geq\, -\frac18\log q .
\end{equation}
In order to derive a complementary upper bound on $\xi_{q,N}$ notice that on
the 
level of events (under the vertical state $\mu_{q,N}^v$), 
\[
\left\{ \widetilde{\omega}\cap J_k^v = \emptyset\right\}\, \supset\, 
\left\{\forall\, \gamma\, \text{exteriour contour of $\widetilde{\omega}$}\ 
J_k^v\cap \text{int} (\gamma ) = \emptyset\right\} .
\]
Indeed, by the definition $J_k^v$ can be intersected only by horizontal 
rods. Let us say that a super-contour $\Gamma$ is intersection 
incompatible with $J_k^v$; $\Gamma\stackrel{i}{\not\sim}J_k^v$, if
 $\Gamma$ contains a contour $\gamma$, such that 
$J_k^v\cap \text{int} (\gamma) \neq \emptyset$. Then, by Lemma~\ref{lem:KP},
\[
\mu_{q,N} \left(\widetilde{\omega}\cap J_k^v = \emptyset \right)
\, \geq\, \exp\left( - \sum_{{\mathcal C}\stackrel{i}{\not\sim}J_k^v}
\abs{\Phi^{\mathrm T} ({\mathcal C})}\right)\, \geq\, e^{-ak} .
\]
Consequently, $\xi_{q,N}\leq a$ and, in view of \eqref{tau_bound} and 
 \eqref{achoice}, the 
proof of Lemma~\ref{hv_gap} is concluded.

\subsection{Sketch of a proof of Theorem~\ref{thm_Domains}}
Consider boxes $V_n^{\underline{k}}$ with periodic boundary conditions.
As before we continue to ignore winding super-contours. In particular the 
notion of exteriour colour is always well defined. Let, therefore, 
$Z_{n,\underline{k}}^{h,\text{per}}$ and 
$Z_{n,\underline{k}}^{v,\text{per}}$ be the partition functions of rod
tilings of $V_n^{\underline{k}}$ with the exteriour colour being fixed
as $h$ (respectively $v$). By \eqref{G_length} and Lemma~\ref{lem:KP},
\begin{equation}
\label{same_periodic}
\left| \log\frac{Z_{n,\underline{k}}^{h,\text{per}}}
{Z_{n,\underline{k}}^{v,\text{per}}}\right|\, \leq\, c_3n^2e^{-c_4 a n} .
\end{equation}
Finally, let $\mu_{n, \underline{k}}^{h, \text{per}}$ and 
$\mu_{n, \underline{k}}^{v, \text{per}}$ be the corresponding Gibbs states. 

The partition functions $Z_{n,\underline{k}}^{h,\text{f}}$ and 
$Z_{n,\underline{k}}^{v,\text{f}}$ of the (exteriour colour) horizontal
and vertical tilings of $V_{n}^{\underline{k}}$ with {\em free} boundary 
conditions are related to $Z_{n,\underline{k}}^{h,\text{per}}$ and 
$Z_{n,\underline{k}}^{v,\text{per}}$ as follows: Set 
\[
J_{n, \underline{k}}^v = \setof{i= (i_1 ,i_2 )\in 
V_n^{\underline{k}}}{i_1=0}\quad \text{and} \quad 
  J_{n, \underline{k}}^h = \setof{i= (i_1 ,i_2 )\in 
V_n^{\underline{k}}}{i_2=0} .
\] 
Then,
\[
\frac{Z_{n,\underline{k}}^{h,\text{f}}}{Z_{n,\underline{k}}^{h,\text{per}}}
\, =\, 
\mu_{n, \underline{k}}^{h, \text{per}}\left( 
\widetilde{\omega}\cap J_{n, \underline{k}}^h=\emptyset\, ;\, 
\widetilde{\omega}\cap J_{n, \underline{k}}^v=\emptyset\right) ,
\]
and, respectively, 
\[
\frac{Z_{n,\underline{k}}^{v,\text{f}}}{Z_{n,\underline{k}}^{v,\text{per}}}
\, =\, \mu_{n, \underline{k}}^{v, \text{per}}\left( 
\widetilde{\omega}\cap J_{n, \underline{k}}^h=\emptyset\, ;\, 
\widetilde{\omega}\cap J_{n, \underline{k}}^v=\emptyset\right) .
\]
By \eqref{tau_qN} and \eqref{xi_qN} the latter probabilities are 
logarithmically asymptotic to 
\[
\exp\left( -n((2k_2 +1)\tau_{q,N} + (2k_1 +1)\xi_{q,N})\right)\quad 
\text{and} \quad
\exp\left( -n((2k_1 +1)\tau_{q,N} + (2k_2 +1)\xi_{q,N})\right) 
\] 
respectively. The claim of Theorem~\ref{thm_Domains} follows now 
from \eqref{same_periodic} and Lemma~\ref{hv_gap}. 

%
%
%

\bibliographystyle{plain}
\bibliography{IVZ06}

\end{document}